\documentclass[12pt,twoside]{article}
\usepackage{amsxtra,amssymb,amsthm,amsmath,latexsym}

\theoremstyle{plain}

\def\oH{\buildrel\circ\over H}
\def\oH1{\buildrel\circ\over H\kern-.02in{}^1}
\def\oH1{\buildrel\circ\over H\kern-.02in{}^1}

\def\d{\delta}
\def\ep{\epsilon}

\def\a{\alpha}

\begin{document}

\title{Two results on ill-posed problems
\thanks{Math subject classification: 47A05, 45A50,35R30}
\thanks{key words:  linear operators, 
ill-posed problems, regularization, 
discrepancy principle }
}

\author{                          
A.G. Ramm\\       
 Mathematics Department, Kansas State University, \\
 Manhattan, KS 66506-2602, USA\\
ramm@math.ksu.edu\\}

\date{}

\maketitle \thispagestyle{empty}

{\bf Abstract} 

Let $A=A^*$ be a linear operator in a Hilbert space $H$.
Assume that equation $Au=f \quad (1)$ is solvable, not
necessarily uniquely, and $y$ is its minimal-norm solution.
Assume that problem (1) is ill-posed. Let $f_\d$,
$||f-f_d||\leq \d$, be noisy data, which are given, while
$f$ is not known.  Variational regularization of problem (1)
leads to an equation $A^*Au+\a u=A^*f_\d$. Operation count
for solving this equation is much higher, than for solving
the equation $(A+ia)u=f_\d \quad (2)$. The first result is
the theorem which says that if $a=a(\d)$, $\lim_{\d \to
0}a(\d)=0$ and $\lim_{\d \to 0}\frac \d {a(\d)}=0$, then the
unique solution $u_\d$ to equation (2), with $a=a(\d),$ has
the property $\lim_{\d \to 0}||u_\d-y||=0$. The second
result is an iterative method for stable calculation of the
values of unbounded operator on elements given with an
error.

1. {\bf Introduction}

 The results of this note are formulated as Theorems 1 and 2
and proved in Sections 1 and 2 respectively. For the notions,
related to ill-posed problems, one may consult \cite{M} and \cite{R}
and the literature cited there.

Let $A=A^*$ be a linear operator in a Hilbert space $H$. Assume that
equation
$$Au=f,
 \eqno{(1)}
$$
is solvable, not necessarily uniquely, and $y$ is its minimal-norm
solution, $y\bot N:=\{u: Au=0\}$. 
Assume that problem (1) is ill-posed. In this case small perturbations of
$f$ may cause large perturbations of the solution to (1) or may
throw $f$ out of the range of $A$. 
Let $f_\d$, $||f-f_d||\leq \d$, be noisy data, which are given, while $f$ 
is not known.
Variational regularization of problem (1) leads to an equation
$A^*Au+\a u=A^*f_\d$, where $A^*=A$ since we assume $A$ to be selfadjoint. 
Operation count for solving this equation
is much higher, than for solving the equation 
$$(A+ia)u=f_\d,
 \eqno{(2)}
$$
The result of this paper is the following theorem.

{\bf Theorem 1.} {\it Let $A=A^*$ be a linear bounded, or densely defined, 
unbounded, selfadjoint operator in a Hilbert space. Assume that 
$a=a(\d)>0$,
$\lim_{\d \to 0}a(\d)=0$ and $\lim_{\d \to 0}\frac \d
{a(\d)}=0$, then the unique solution $u_\d$ to equation (2)
with $a=a(\d)$ has the property 
$$\lim_{\d \to 0}||u_\d-y||=0.
 \eqno{(3)}
$$
} 
\vskip .07in
Why should one be interested in the above theorem? 
The answer is: because the solution to equation (2)
requires less operations than the solution of the  
equation $(A^*A+\a I)u=A^*f_\d$ basic for the variational 
regularization method for stable solution of equation (1).
Here $I$ is the identity operator. 
Also, a discretized version of (2) leads to matrices whose condition 
number is of the order of square root of the condition number
of the matrix corresponding to the operator $A^*A+\a(\d)I$, where 
$I$ is the identity operator.

{\bf Proof of Theorem 1.}  One has
$$||u_{a,\d}-y||\leq ||(A+ia)^{-1}(f_\d -f)||+||(A+ia)^{-1}Ay-y||\leq 
\frac \d a+a||(A+ia)^{-1}y||.
\eqno{(4)}
$$
Moreover, 
$$\lim_{a \to 0}a^2||(A+ia)^{-1}Ay-y||^2=\lim_{a \to 
0}a^2\int_{-\infty}^{\infty} \frac{d(E_sy,y)}{s^2+a^2}=0,
\eqno{(5)}
$$
where we have used the spectral theorem, $E_s$ is the 
resolution of the identity corresponding to the selfadjoint operator $A$,
and we have taken into account that 
$$\lim_{a \to 0}a^2\int_{-0}^{0} \frac{d(E_sy,y)}{s^2+a^2}=0$$ 
because 
$y\bot N$.
From formulas (4) and (5) one concludes that if
 $a=a(\d)>0$,
$\lim_{\d \to 0}a(\d)=0$ and $\lim_{\d \to 0}\frac \d{a(\d)}=0$, 
then the unique solution $u_\d$ to equation (2)
with $a=a(\d)$ satisfies equation (3).

Theorem 1 is proved. \hfill $\Box$

2. {\bf Calculation of values of unbounded operators}

Assume that $A$ is a densely defined closed 
linear operator in $H$. We do not assume in this Section that $A$ is 
selfadjoint. If $f\in D(A)$, then we want to compute $Af$ given noisy
data $f_\d$, $||f_\d-f||\leq \d$. Note that $f_\d$ may not belong to 
$D(A)$. The problem of stable calculation of $Af$ given the data
$\{f_\d, \d, A\}$ is ill-posed. It was studied in the literature
(see, e.g., \cite{M}) by a variational regularization method.
Our aim is to reduce this problem to a standard equation with a 
selfadjoint bounded operator $0\leq B\leq I$,  and solve this equation 
stably by an iterative method.
 
Let $v=Af$. This relation is equivalent to 
$$Bv=Ff,
\eqno{(6)}
$$
where $B:=(I+Q)^{-1}$, $F:=BA$, $Q=AA^*$ is a densely defined, 
non-negative, selfadjoint operator,
the range of $I+Q$ is the whole space $H$, and $B$ is a selfadjoint 
operator, $0\leq B\leq I$, where the inequalities are understood in the 
sense of quadratic forms, e.g., $B\geq 0$ means $(Bg,g)\geq 0$ for all
$g\in H$, and $F:=(I+Q)^{-1}A$.

{\bf Lemma 1.}  (see \cite{R1}) {\it The operator $(I+Q)^{-1}A$, 
originally defined on $D(A)$, is closable. Its closure is a bounded, 
defined on all of $H$ linear operator with the norm $\leq \frac 12$.
One has $(I+Q)^{-1}A=A(I+T)^{-1}$, where $T=A^*A$ is a non-negative, 
densely defined selfadjoint operator, and $||A(I+T)^{-1}||\leq \frac 1 
2.$}

If $f_\d$ is given in place of $f$, then we stably solve equation (6) for 
$v$ using the following iterative process:
$$v_{n+1}=(I-B)v_n+Ff_\d, \quad v_0\bot N^*,
\eqno{(7)}
$$
where $ N^*:=\{u: A^*u=0\}$. Let $y$ be the unique minimal-norm solution 
to equation (6), $By=Ff$. Note that $y=Hy+Ff$, where $H:=I-B$.

{\bf Theorem 2.} {\it If $n=n(\d)$ is an integer, $\lim_{\d\to 
0}n(\d)=\infty$ and
$\lim_{\d\to 0}[\d n(\d)]=0$, then 
$$\lim_{\d\to 0}||v_\d-y||=0,
\eqno{(8)}
$$
where $v_\d:=v_{n(\d)}$, and $v_n$ is defined in (7).}

{\bf Proof of Theorem 2.} From (7) one gets $v_{n+1}=\sum_{j=0}^n H^jFf_\d 
+H^{n+1}u_0$, where $H:=I-B$. One has $y=Hy+Ff$. Let $w_n:=v_n-y$.
Then $w_{n}=\sum_{j=0}^{n-1} H^jFg_\d +H^{n}w_0$, where $g_\d:=f_\d-f$,
and $w_0$ is an arbitrary element such that $w_0\bot N^*$.
Since $0\leq H \leq I$,$||F||\leq \frac 1 2$, and $||g_\d||\leq \d$, one 
gets
$$ ||w_n||\leq \frac {n\d}2 +[\int_0^1(1-s)^{2n}d(E_s,w_0,w_0)]^{1/2},
\eqno{(9)}
$$
where $E_s$ is the resolution of the identity corresponding to
the selfadjoint operator $B$.

If $w_0\bot N^*$, then 
$$ 
\lim_{h\to 1}\int_h^1(1-s)^{2n}d(E_s,w_0,w_0)=||Pw_0||^2,
\eqno{(10)}
$$
where $P$ is the orthoprojector onto the subspace $\{u: Bu=u\}=\{u: 
Qu=0\}=N^*$, and $Pw_0=0$ because $w_0\bot N^*$ by the assumption.
The conclusion of Theorem 2 can now be derived. Given an arbitrary small 
$\ep>0$,
find $h$ sufficiently close to $1$ such that 
$\int_h^1(1-s)^{2n}d(E_s,w_0,w_0)<\ep$. Fix this $h$ and find 
$n=n(\d)$, sufficiently large, so that $\d n(\d)<\ep$ and, at the same 
time,
$(1-h)^{2n(\d)}<\ep$. This is possible if $\d$ is sufficiently small,
because $\lim_{\d\to 0}n(\d)=\infty$ and
$\lim_{\d\to 0}[\d n(\d)]=0$. Then $\int_0^1(1-s)^{2n}d(E_s,w_0,w_0)\leq
\int_0^h(1-s)^{2n}d(E_s,w_0,w_0)+\int_h^1(1-s)^{2n}d(E_s,w_0,w_0)<\ep$,
and inequality (9) shows that (8) holds. Theorem 2 is proved. $\hfill 
\Box.$

{\bf Remark 1.} It is not possible to estimate the rate of convergence in 
(8) without making additional assumptions on $y$ or on $f$.
In \cite{R} one can find examples illustrating similar statements 
concerning various methods for solving ill-posed problems.


\begin{thebibliography}{10}

\bibitem{M} V. Morozov, Methods of solving incorrectly posed problems,
Springer Verlag, New York, 1984.
 
\bibitem{R} A. G. Ramm,  Inverse Problems, Springer, 
New York, 2005.

\bibitem{R1} A. G. Ramm, On unbounded operators and applications,
(submitted)

\end{thebibliography}
\end{document}